\documentclass[a4paper, 12pt]{extarticle}
\usepackage{epsfig}
\usepackage{graphicx}
\usepackage{psfrag}
\usepackage{ifpdf}
\usepackage{cite}
\usepackage{siunitx}
\usepackage[utf8x]{inputenc}
\usepackage[top=1in, bottom=1in]{geometry}
\usepackage{fancyhdr}
\usepackage{lastpage}
\usepackage{oldstyle}
\usepackage{textcomp}
\usepackage{color}
\usepackage{placeins}
\usepackage{float}
\usepackage{tabularx,colortbl}
\graphicspath{{figures/Appendix/}}
\usepackage{amssymb} 
\usepackage{amsmath} 
\usepackage{subfig}

\def \Cfi {C_{f}}
\def \Csi {C_{s}}

\DeclareMathOperator*{\argmin}{arg\,min}

\begin{document}

\begin{center}
 \LARGE Projection Onto Convex Sets (POCS) Based Signal Reconstruction Framework with an associated cost function\\[10pt]
 \normalsize\today
\end{center}

\begin{center}
\normalsize Mohammad Tofighi, Kivanc Kose$^*$, A. Enis Cetin\\
Dept. of Electrical and Electronic Engineering, Bilkent University,  Ankara, Turkey\\
$^*$Dermatology Department, Memorial Sloan Kettering Cancer Center, New York, USA\\
tofighi@ee.bilkent.edu.tr, $^*$kosek@mskcc.org, cetin@bilkent.edu.tr\\
\end{center}

\begin{abstract}
A new signal processing framework based on the projections onto convex sets (POCS) is developed for solving convex optimization problems. The dimension of the minimization problem is lifted by one and the convex sets corresponding to the epigraph of the cost function are defined. If the cost function is a convex function in $R^N$ the corresponding epigraph set is also a convex set in $R^{N+1}$. The iterative optimization approach starts with an arbitrary initial estimate in $R^{N+1}$ and orthogonal projections are performed onto epigraph set in a sequential manner at each step of the optimization problem. The method provides globally optimal solutions in total-variation (TV), filtered variation (FV), $\ell_1$, $\ell_1$, and entropic cost functions. New denoising and compressive sensing algorithms using the TV cost function are developed. The new algorithms do not require any of the regularization parameter adjustment. Simulation examples are presented.\end{abstract}

\section{Introduction}
\label{sec:Introduction}
A de-noising method based on a new POCS framework is introduced. In standard POCS approach, only a common point of convex constraint sets is determined. It is shown that it is possible to use a convex cost function in this framework \cite{Rud92, Bar07,Can08, Kos12,Gunay}.

Bregman developed iterative methods based on the so-called Bregman distance  to solve convex optimization problems \cite{Bre67}. In Bregman's approach, it is necessary to perform a D-projection (or Bregman projection) at each step of the algorithm and it may not be easy to compute the Bregman distance in general \cite{Yin08,Kiv12,Gunay}.

In this article Bregman's older projections onto convex sets (POCS) framework \cite{Bregman,You82} is used to solve convex optimization problems instead of the Bregman distance approach. In the ordinary POCS approach the goal is simply to find a vector which is in the intersection of convex sets \cite{Kose2013,You82,Her95,Cen12,Sla08,Cet03,Cetin94,Cetin89,Kose11,Cen81,Sla09,The11,censor1987optimization,Tru85,Com04,Com93,Kim92,yamada2011minimizing,censor1987some,Sez82,censor1992proximal,Tuy81,censor1981row,censor1991optimization,Ros13,Cevher2}. In each step of the iterative algorithm an orthogonal projection is performed onto one of the convex sets. Bregman showed that successive orthogonal projections converge to a vector which is in the intersection of all the convex sets. If the sets do not intersect iterates oscillate between members of the sets \cite{Gub67,Com12,Cet97}. Since  there is no need to compute the Bregman distance in standard POCS, it found applications to many practical problems.

In our approach the dimension of the signal reconstruction or restoration problem is lifted by one and  sets corresponding to a given convex cost function are defined. This approach is graphically illustrated in Fig. \ref{app:convex}. If the cost function is a convex function in $R^N$ the corresponding epigraph set is also a convex set in $R^{N+1}$. As a result the convex minimization problem is reduced to finding a specific member (the optimal solution) of the set corresponding to the cost function. As in ordinary POCS approach the new iterative optimization method starts with an arbitrary initial estimate in $R^{N+1}$ and an orthogonal projection is performed onto one of the sets. After this vector is calculated it is projected onto the other set. This process is continued in a sequential manner at each step of the optimization problem.
The method provides globally optimal solutions in total-variation, filtered variation, $\ell_1$, and entropic function based cost functions because they are convex cost functions.

The article is organized as follows. In Section \ref{sec:Convex Minimization}, the convex minimization method based on the POCS approach is introduced. In Section \ref{sec:Denoising using POCS}, a new denoising method based on the convex minimization approach introduced in Section \ref{sec:Convex Minimization}, is presented.This new approach uses supporting hyperplanes of the TV function and it does not require a regularization parameter as in other TV based methods. Since it is very easy to perform an orthogonal projection onto a hyperplane this method is computationally implementable for many cost functions without solving any nonlinear equations. In Section~\ref{sec:Simulation Results}, we present the simulation results and some denoising examples.

\section{Convex Minimization Using The Epigraph set}
\label{sec:Convex Minimization}
Let us first consider a convex minimization problem
\begin{equation}
\label{app:eq:c5}
\underset{\mathbf{w}\in \mathbb{R}^N}{\text{min}} f(\mathbf{w}),
\end{equation}
where $f:\mathbb{R}^N \rightarrow \mathbb{R}$ is a convex function.
We increase the dimension by one to define the following sets in $\mathbb{R}^{N+1}$ corresponding to the cost function $f(\mathbf{w})$ as follows:
\begin{equation}
\label{app:eq:c6}
\text{C}_f = \{\underline{\mathbf{w}} = [\mathbf{w}^T~y]^T : \mathrm{~} y\geq f(\mathbf{w})\},
\end{equation}
which is the set of $N+1$ dimensional vectors whose $(N+1)^{st}$ component $y$ is greater than $f(\mathbf{w})$. This set $C_{f}$ is called the epigraph of $f$. We use bold face letters for $N$ dimensional vectors and underlined bold face letters for $N+1$ dimensional vectors, respectively.

\begin{figure}[ht!]
\begin{center}
\noindent
\includegraphics[width=80mm]{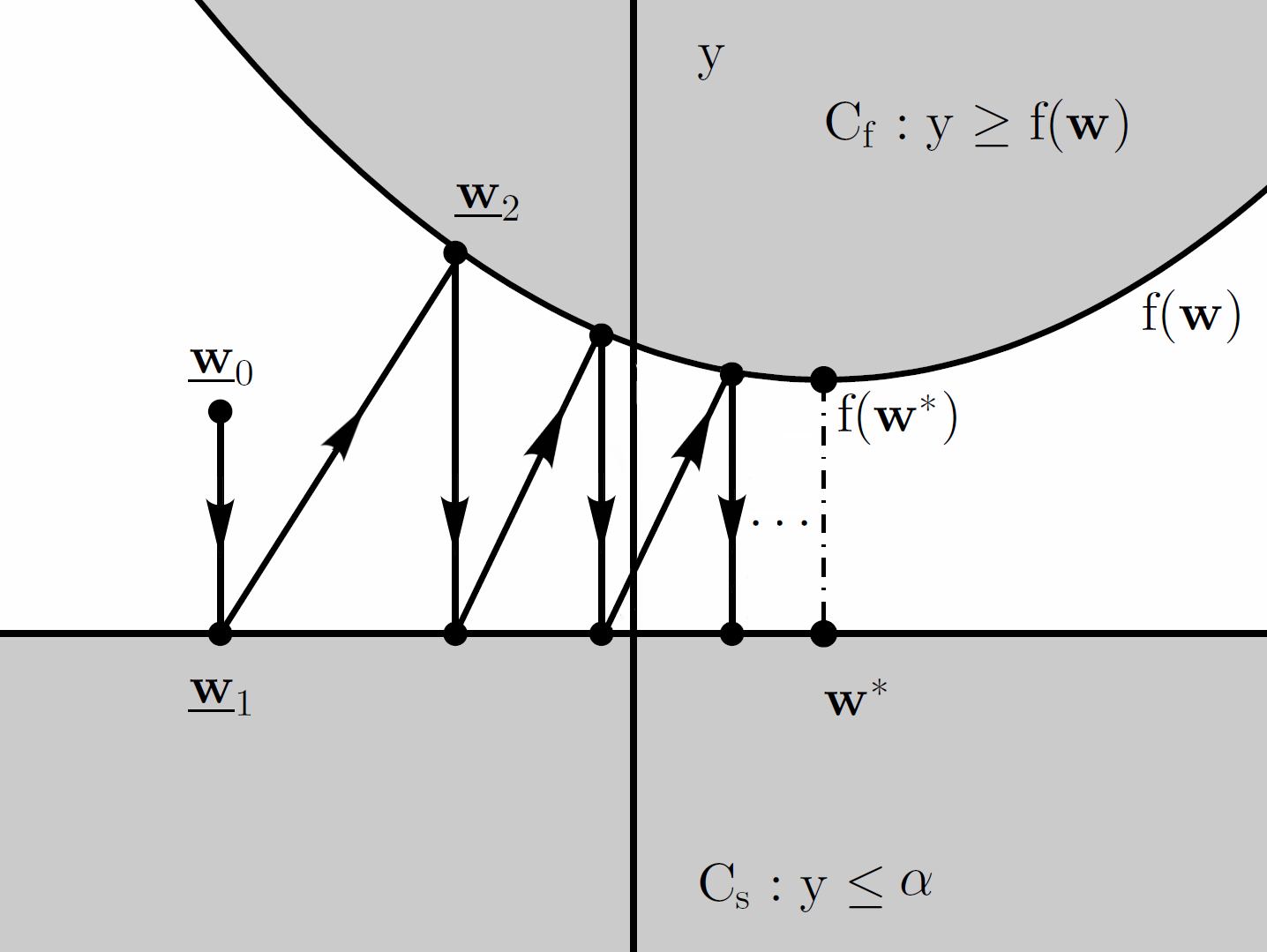}
\caption[Two projecting convex sets.]{Two convex sets $\Cfi$ and $\Csi$ corresponding to the cost function $f$. We sequentially project an initial vector $\underline{\mathbf{w}}_0$ onto $\Csi$ and C$\Cfi$ to find the global minimum which is located at $\mathbf{w}^*$.}
\label{app:convex}
\end{center}
\end{figure}

The second set that is related with the cost function $f(\mathbf{w})$ is the level set:
\begin{equation}
\label{app:eq:c7}
\text{C}_s = \{\underline{\mathbf{w}} = [\mathbf{w}^T~y]^T : ~ y\leq \alpha , ~\underline{\mathbf{w}} \in \mathbb{R}^{N+1}\},
\end{equation}
where $\alpha$ is a real number. Here it is assumed that $f(\mathbf{w})\geq \alpha$ for all $f(\mathbf{w})\in \mathbb{R}$ such that the sets $\Cfi$ and $\Csi$ do not intersect. They are both closed and convex sets in  $\mathbb{R}^{N+1}$. Sets $\Cfi$ and $\Csi$ are graphically illustrated in Fig. \ref{app:convex} in which $\alpha=0.$

The POCS based minimization algorithm starts with an arbitrary $\underline{\mathbf{w}}_0 =[ \mathbf{w}_0^T ~ y_0]^T \in \mathbb{R}^{N+1}$. We project $\underline{\mathbf{w}}_0$ onto the set $\Csi$ to obtain the first iterate $\underline{\mathbf{w}}_1$ which will be,
\begin{equation}
\underline{\mathbf{w}}_1 = [\mathbf{w}_0^T~0]^T,
\end{equation}
where $\alpha=0$ is assumed as in Fig. \ref{app:convex}. Then we project $\underline{\mathbf{w}}_1$ onto the set $\Cfi$. The new iterate  $\underline{\mathbf{w}}_2$ is determined by minimizing the distance between $\underline{\mathbf{w}}_1$ and $\Cfi$, i.e.,

\begin{equation}
\label{app:eq:convex}
\underline{\mathbf{w}}_2 = \text{arg} \underset{\underline{\mathbf{w}}\in \text{C}_{\mathrm{s}}}{\text{min}} \|\underline{\mathbf{w}}_1 - \underline{\mathbf{w}}\|.
\end{equation}

Equation \ref{app:eq:convex} is the ordinary orthogonal projection operation onto the set $\mathrm{C_f} \in \mathbb{R}^{N+1}$.
To solve the problem in Eq. \ref{app:eq:convex} we do not need to compute the Bregman's so-called D-projection. After finding $\underline{\mathbf{w}}_2$, we perform the next projection onto the set $\Csi$ and obtain $\underline{\mathbf{w}}_3$ etc. Eventually iterates oscillate between two nearest vectors of the two sets $\Csi$ and $\Cfi$. As a result we obtain
\begin{equation}
\label{app:eq:convex2}
\underset{n \rightarrow \infty}{\text{lim}} \underline{\mathbf{w}}_{2n} = [\mathbf{w}^*~f(\mathbf{w}^*)]^T,
\end{equation}
where $\mathbf{w}^*$ is the N dimensional vector minimizing $f(\mathbf{w})$. The proof of Eq. \ref{app:eq:convex2} follows from Bregman's POCS theorem \cite{Bregman,Gub67}. It was generalized to non-intersection case by Gubin et. al \cite{Gub67,Cen12},\cite{Com12}. Since the two closed and convex sets $\Csi$ and $\Cfi$ are closest to each other at the optimal solution case, iterations oscillate between the vectors $[\mathbf{w}^*~f(\mathbf{w}^*)]^T$ and $[\mathbf{w}^*~0]^T$ in $R^{N+1}$ as $n$ tends to infinity. It is possible to increase the speed of convergence by non-orthogonal projections \cite{Com93}.

If the cost function $f$ is not convex and have more than one local minimum then the corresponding set $\Cfi$ is not convex in $R^{N+1}$. In this case iterates may converge to one of the local minima.

Consider the standard LASSO based denoising \cite{Cevher1}:
\begin{equation}
\label{app:eq:LASSO}
{{\text{min}} \frac{1}{2}\|{\mathbf{v}} - {\mathbf{w}}\|}^2_2 + \lambda \|\textbf{w}\|_1,
\end{equation}
where $\textbf{v}$ is the corrupted version of \textbf{w}. Since the cost function
\begin{equation}
\label{app:eq:costLASSO}
f(\textbf{w}) = \frac{1}{2}\|{\mathbf{y}} - {\mathbf{w}}\|^2_2 + \lambda \|\textbf{w}\|_1,
\end{equation}
is a convex function, the framework introduced in this section can solve this problem in an iterative manner. One weakness of this approach is that the smoothing or regularization parameter $\lambda$ has to be specified or manually selected.

\section{De-noising Using POCS}
\label{sec:Denoising using POCS}
\begin{figure}[ht!]
\begin{center}
\noindent
\includegraphics[width=90mm]{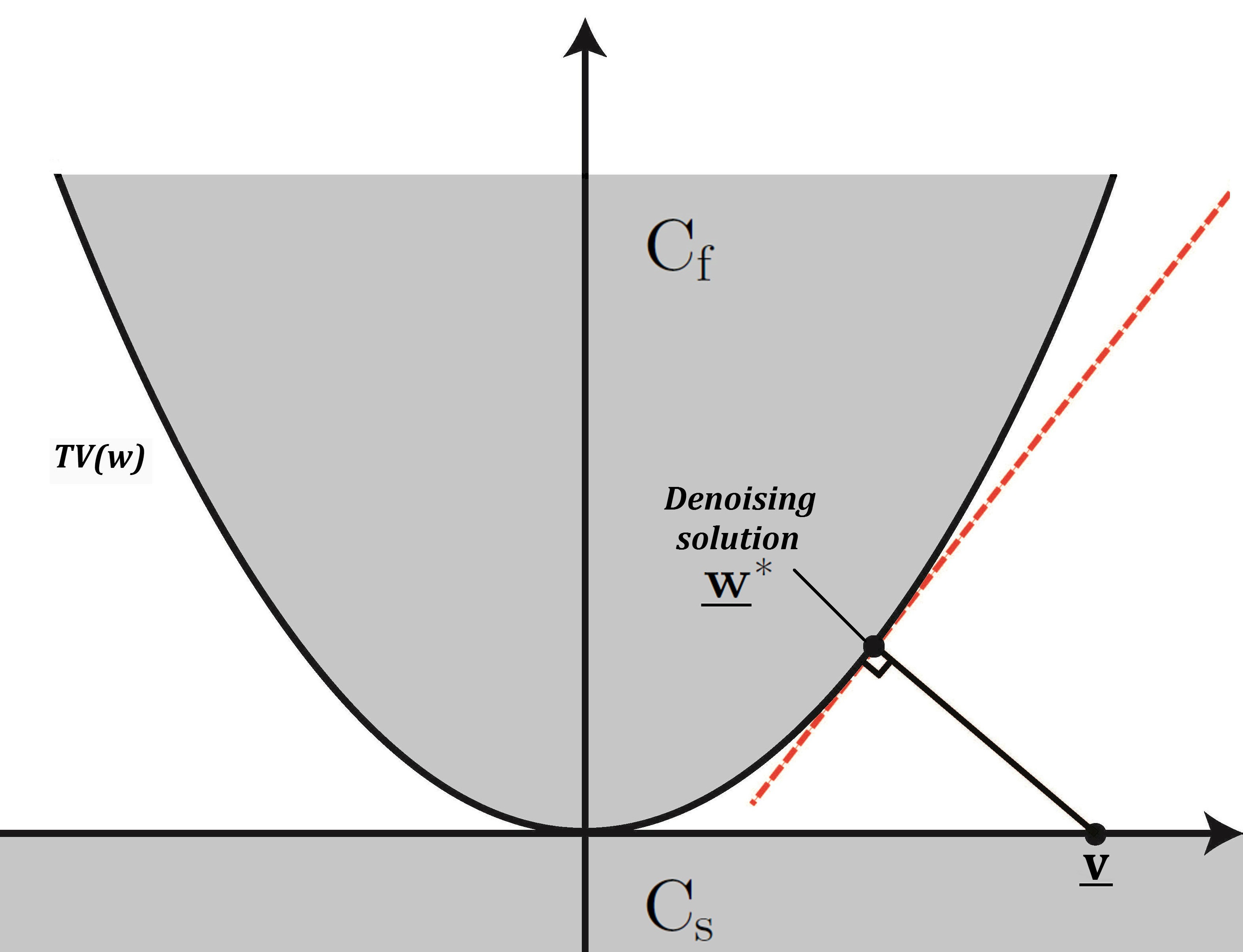}
\caption{Graphical representation of the minimization of (\ref{app:eq:6}). The corrupted observation vector $\textbf{v}$ is projected onto the set $\Cfi$. TV(\textbf{w}) is zero for $\textbf{w}=[0, 0,...,0]^{T}$ or when it is a constant vector.}
\label{app:convexSol}
\end{center}
\end{figure}

In this section, we present a new method of denoising, based on TV and FV. Let the noisy signal be \textbf{v}, and the original signal or image be $\textbf{w}_{0}$. Suppose that the observation model is the additive noise model:
\begin{equation}
\textbf{v} = \textbf{w}_{0} + \boldsymbol{\eta},
\end{equation}
where $\boldsymbol{\eta}$ is the additive noise. In this approach we solve the following problem for denoising:
\begin{equation}
\label{app:eq:6}
\underline{\mathbf{w}}^{\star} = \text{arg} \underset{\underline{\mathbf{w}}\in \text{C}_{\mathrm{f}}}{\text{min}} \|\underline{\mathbf{v}} - \underline{\mathbf{w}}\|^{2},
\end{equation}
where, $\underline{\mathbf{v}}$ = [$v^{T}$\ 0] and $\Cfi$ is the epigraph set of TV or FV in $R^{N+1}$. The TV function which we used for discrete image $\textbf{w} = [w^{i,j}]~~0\leq i,j\leq M-1~\in f:\mathbb{R}^{M\times M}$ is as follows:
\begin{equation}
\label{app:eq:7}
TV(\textbf{w}) = \sum_{i,j=1}^{M} |w^{i+1, j} - w^{i, j}| + \sum_{i,j=1}^{M} |w^{i, j+1} - w^{i, j}|.
\end{equation}
The minimization problem is essentially the orthogonal projection onto the set $C_{f,i}$. This means that we select the nearest vector $\underline{\mathbf{w}}^{\star}$ on the set $C_{f,i}$ to \textbf{v}. This is graphically illustrated in Fig. \ref{app:convexSol}.

Actually, Combettes and Pesquet and other researchers including us used a similar convex set in denoising and other signal restoration applications [].  The following convex set describes all signals whose TV is bounded by $\epsilon$:
\begin{equation}
\label{convex_set}
C_b = \{  {\bf w}: TV( {\bf w} ) <= \epsilon \}
\end{equation}
The parameter $\epsilon$ is a fixed upper bound on the total variation of the signal and it has to be determined in an ad hoc manner a priori. On the other hand we do not specify a prescribed number on the TV in the epigraph approach. The upperbound on TV is automatically determined by the orthogonal projection
according to the location of the corrupted signal as shown in Fig. \ref{app:convexSol}.

The denoising solution $\underline{\mathbf{w}}^{\star}$ has the lowest total variation on the line [$\underline{\mathbf{v}}$, $\underline{\mathbf{w}}^{\star}$]. In current TV based denoising methods \cite{Chambolle, Com11} the following cost function is used:

\begin{equation}
\label{app:eq:cost}
{{\text{min}} \|\mathbf{v}} - {\mathbf{w}}\|^2_2 + \lambda \text{TV}(\textbf{w}).
\end{equation}
The solution of this problem can be obtained using the method that we discussed in Section \ref{sec:Convex Minimization}. One problem with this approach is the estimation of the regularization parameter $\lambda$. One has to determine the $\lambda$ in an ad hoc manner or by visual inspection. On the other hand we do not require any parameter adjustment in (\ref{app:eq:6}).

\begin{figure}[ht!]
\begin{center}
\noindent
\includegraphics[width=90mm]{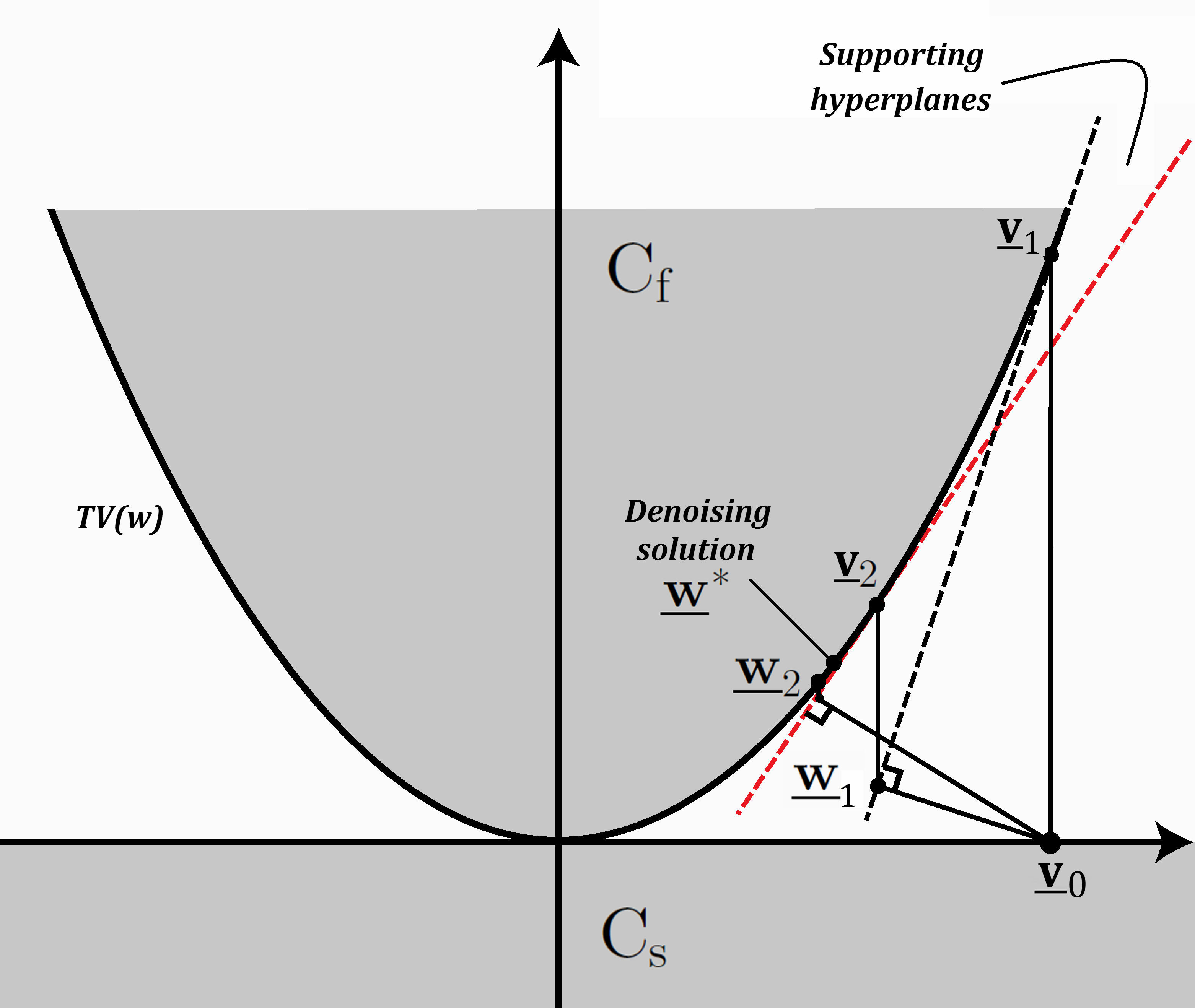}
\caption{Graphical representation of the minimization of (\ref{app:eq:6}), using projection onto the supporting hyperplanes of $\Cfi$.}
\label{app:convex1}
\end{center}
\end{figure}

The denoising solution in (\ref{app:eq:6}) can be found by performing successive orthogonal projection onto supporting hyperplanes of the epigraph set $\Cfi$. In the first step we calculated TV($\mathbf{v}$). We also calculate the surface normal at $\underline{\mathbf{v}}$ = [$\mathbf{v}^{T}$ \ TV($\mathbf{v}$)] in $R^{N+1}$ and determine the equation of the supporting hyperplane at [$\mathbf{v}^{T}$ \ TV($\mathbf{v}$)].
We project $\underline{\mathbf{v}}$ = [$\mathbf{v}^{T}$ \ 0] onto this hyperplane and obtain $\underline{\mathbf{w}}_1$ as our first estimate as shown in Fig. \ref{app:convex1}. In the second step we project $\underline{\mathbf{w}}_1$ onto the set $\Csi$ by simply making its last component zero. We calculate the TV of this vector and the surface normal, and the supporting hyperplane as in the previous step. We project $\underline{\mathbf{v}}$ onto the new supporting hyperplane, etc.

\begin{figure}[ht!]
\begin{center}
\noindent
\includegraphics[width=90mm]{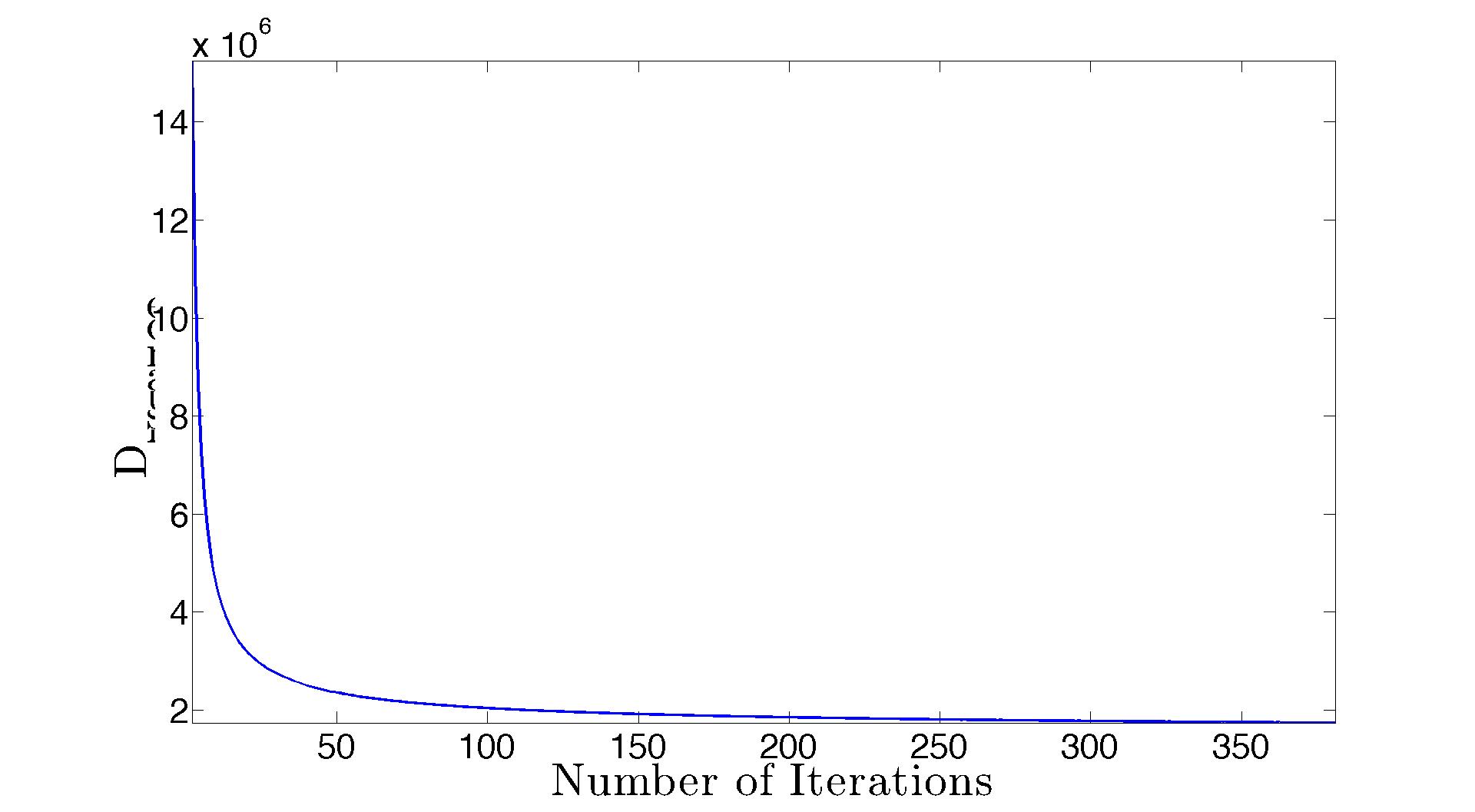}
\caption{Orthogonal distance from v to the epigraph of TV in each iteration.}
\label{app:graphdist}
\end{center}
\end{figure}

The sequence of iterations obtained in this manner converges to a vector in the intersection of $\Csi$ and $\Cfi$. In this problem the sets $\Csi$ and $\Cfi$ intersect because $TV(\textbf{w})=0$ for $\textbf{w}=[0, 0,...,0]^{T}$ or for a constant vector. However, we do not want to find a trivial constant vector in the intersection of $\Csi$ and $\Cfi$. We calculate the distance between $\underline{\mathbf{v}}$ and $\underline{\mathbf{w}}_i$ at each step of the iterative algorithm described in the previous paragraph. This distance ${ \|\underline{\mathbf{v}} - \underline{\mathbf{w}}_i\|}^2_2$ initially decreases and starts increasing as $i$ increases. Once we detect the increase we perform some refinement projections to obtain the solution of the de-noising problem. A typical convergence graph is shown in Fig. \ref{app:graphdist} for the ``note" image.

\section{Compressive Sensing}
The most common method used in compression applications is transform coding. The signal $\textbf{x}$ is transformed into another domain defined by the transformation matrix $\boldsymbol{\psi}$ . The transformation procedure is simply finding the inner product of the signal $\textbf{x}$ with the rows $\psi_{i}$ of the transformation matrix $\boldsymbol{\psi}$ represented as follows:
\begin{equation}
\label{app:eq:eq1}
s_{l} = \langle \textbf{x}, \psi_{l} \rangle \quad i = 1, 2, ..., N,
\end{equation}
where \textbf{x} is a column vector of size $N$.

The discrete time signal can be reconstructed from its transform coefficients $s_{l}$ as follows:
\begin{equation}
\label{app:eq:eq2}
\textbf{x} = \sum_{l=1}^{N} s_l.\psi_{l}\quad or\quad \textbf{x} = \psi.s,
\end{equation}
where \textbf{s} is a vector containing the transform domain coefficients $s_{l}$.

The basic idea in digital waveform coding is that the signal should be approximately reconstructed from only a few of its non-zero transform coefficients. In most cases, including the JPEG image coding standard, the transform matrix $\psi$ is chosen in such a way that the new signal s is efficiently represented in the transform domain with a small number of coefficients. A signal $x$ is compressible, if it has only a few large amplitude $s_{l}$ coefficients in the transform domain and the rest of the coefficients are either zeros or negligibly small-valued.

In a compressive sensing framework, the signal is assumed to be K-Sparse in a transformation domain, such as the wavelet domain or the DCT (Discrete Cosine Transform) domain. A signal with length $N$ is K-Sparse if it has at most $K$ non-zero and $(N − K)$ zero coefficients in a transform domain. The case of interest in Compressive Sensing (CS) problems is when $K\ll N$, i.e., sparse in the transform domain.

The CS theory introduced in \cite{Candes1, Bar07, candes2006compressive, Candes2, Donoho1} provides answers to the question of reconstructing a signal from its compressed measurement vector \textbf{v}, which is defined as follows:
\begin{equation}
\label{app:eq:eq3}
\textbf{v} = \psi x = \psi.\phi.s = \theta.s,
\end{equation}
where $\phi$ is the $M\times N$ measurement matrix and $M\ll N$. The reconstruction of the original signal x from its compressed measurements y cannot be achieved by simple matrix inversion or inverse transformation techniques. A sparse solution can be obtained by solving the following optimization problem:
\begin{equation}
\label{app:eq:eq4}
s_{p} = \argmin \|s\|_{0}\quad such~that\quad \theta.s = y.
\end{equation}

However, this problem is an NP-complete optimization problem; therefore, its solution can not be found easily. It is also shown in \cite{Candes1, Bar07, candes2006compressive, Candes2} that it is possible to construct the $\phi$ matrix from random numbers, which are i.i.d. Gaussian random variables. In this case, the number of measurements should be chosen as $cK log(\frac{N}{K}) < M \ll N$ \cite{Candes1}, \cite{Bar07}. With this choice of the measurement matrix, the optimization problem (\ref{app:eq:eq4}) can be approximated by $\ell_1$ norm minimization as:
\begin{equation}
\label{app:eq:eq5}
s_{p} = \argmin \|s\|_{1}\quad such~that\quad \theta.s = y.
\end{equation}

Instead of solving the original CS problem in (\ref{app:eq:eq4}) or (\ref{app:eq:eq5}), several researchers developed methods to reformulate those and approximate the solution through these new formulations. For example, in \cite{Ji1}, the authors developed a Bayesian framework and solved the CS problem using Relevance Vector Machines (RVM). Some researchers replaced the objective function of the CS optimization in (\ref{app:eq:eq4}), (\ref{app:eq:eq5}) with a new objective function to solve the sparse signal reconstruction problem \cite{Chartrand1, Rezaii}. One popular approach is replacing $\ell_0$ norm with $\ell_p$ norm, where p $\in$ (0, 1) or even with the mix of two different norms as in  \cite{Chartrand1, Rezaii,Chartrand2, Kowalski}. However, in these cases, the resulting optimization problems are not convex. Several studies in the literature addressed $\ell_p$ norm based non-convex optimization problems and applied their results to the sparse signal reconstruction example \cite{ehler2011shrinkage, Achim, tzagkarakis}.
We use the epigraph of $\ell_1$ norm cost function and the TV cost function  together with the measurement hyperplanes to solve this problem. Let
\begin{equation}
\label{app:eq:6}
C_f =  \{  f(w) \leq y \}
\end{equation}
where $f(w)$ is the cost function representing the $\ell_1$ norm or the TV function, and the measurement hyperplane sets are defined as follows:
\begin{equation}
C_i =  \{ w.\phi_{i} = v_{i} \} , \quad i=1,2,...,L
\end{equation}
All of the above sets are closed and convex sets. Therefore it is possible to device an iterative signal reconstruction algorithm in $R^{N+1}$ by making successive orthogonal projections onto the sets $C_s$  and $C_i$.

In this case we replace the level set $C_s$ in Section $\ref{sec:Convex Minimization}$ with the measurement hyperplanes. Since the hyperplanes form an undetermined set of equations their intersection $C_{int} = \cap_i C_i $ is highly unlikely to be an empty set but the intersection of hyperplanes $C_{int}$ may not intersect with the epigraph set $C_f$. This scenario has not been studied in POCS theory to the best of our knowledge [] but it is very similar to the scenario that we discussed in Section $\ref{sec:Convex Minimization}$.

We solve this problem in an iterative manner by performing successive orthogonal projection onto hyperplane corresponding to measurements  $v_{i} = w.\phi_{i}$ for $i = 1, 2, ..., N$, followed by an orthogonal projection onto the epigraph set. As we pointed out in the previous paragraph this case has not been studied in the literature. But it is intuitively clear that iterates oscillate between a vector  in the intersection of hyperplanes $C_{int}$ and the closest vector of the epigraph set $C_f$ to the  intersection similar to the approach introduced in Section  $\ref{sec:Convex Minimization}$. Therefore we essentially obtain a solution to the following problem:
\begin{equation}
\label{minimiz}
min || w_f - w_{int} ||_2
\end{equation}
where $w_f\in C_f $ and $w_{int}\in C_{int}$.

If the sets $C_f$ and $C_{int}$ intersect the iterates converge to a vector in the $C_{int}\cap C_f$ by Bregman's POCS theorem.

In our approach its also possible to define a smoothing parameter in both denoising and compressive sensing  solutions as well. The epigraph set $C_f$ can be modified as follows:
\begin{equation}
\label{app:eq:c6}
\text{C}_{f, \alpha} = \{\underline{\mathbf{w}}: \mathrm{~} y\geq \alpha TV(\mathbf{w})\}.
\end{equation}
The choice of the parameter $\alpha > 1$ provides smoother solution than usual and $\alpha < 1$ relaxes the smoothing constraint. Its experimentally observed that $\alpha = 1$ usually provides better de-noising results than optimally selected $\lambda$ values in standard TV denoising in \cite{Chambolle}. Simulation examples are presented in the next section.

\section{Simulation Results}
\label{sec:Simulation Results}
We present de-noising examples in Section \ref{subsec:Simulation De-noising} and compressive sensing examples in Sensing \ref{subsec:Simulation Compressive sensing}.

\subsection{De-noising}
\label{subsec:Simulation De-noising}
\begin{figure}[ht]
\centering
\subfloat[]
{\label{fig:note_org}\includegraphics[width=70mm]{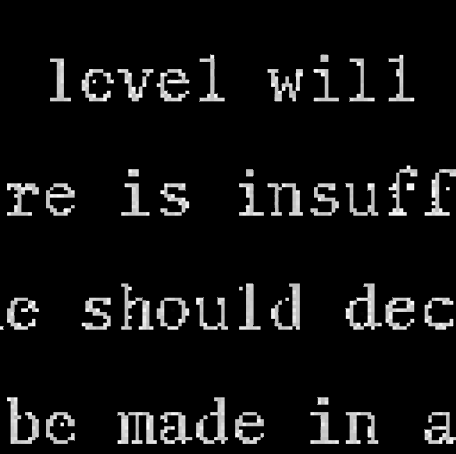}} \quad
\subfloat[]
{\label{fig:note_noise}\includegraphics[width=70mm]{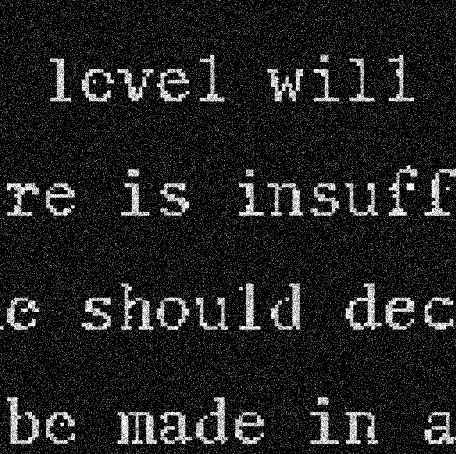}} \qquad
\subfloat[]
{\label{fig:note_denoised_POCS}\includegraphics[width=70mm]{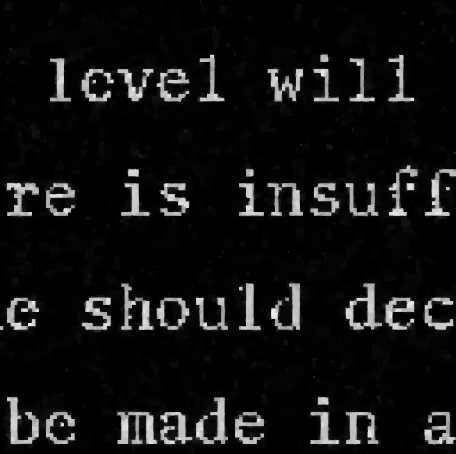}} \quad
\subfloat[]
{\label{fig:note_denoised_chombolle}\includegraphics[width=70mm]{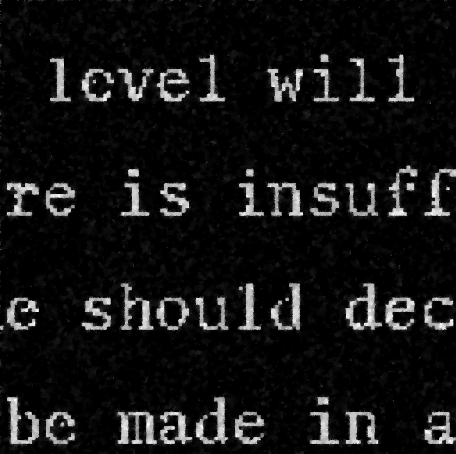}} \qquad
\caption{Sample images used in our experiments (a) Original ``Note" image, (b) ``Note" image corrupted with Gaussian noise with $\delta = 45$, (c) Denoised ``Note" image, using POCS algorithm; SNR = 15.08 dB, (d) Denoised ``Note" image, using Chambolle's algorithm; SNR = 12.78 dB.}
\label{fig:note_example}
\end{figure}

\begin{figure}[ht]
\centering
\subfloat[]
{\label{fig:flower_org}\includegraphics[width=70mm]{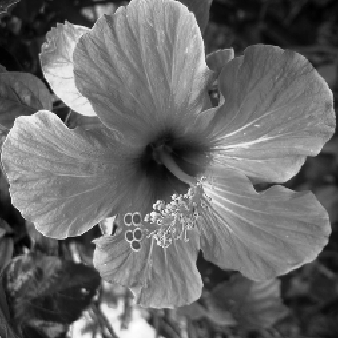}} \quad
\subfloat[]
{\label{fig:flower_noise}\includegraphics[width=70mm]{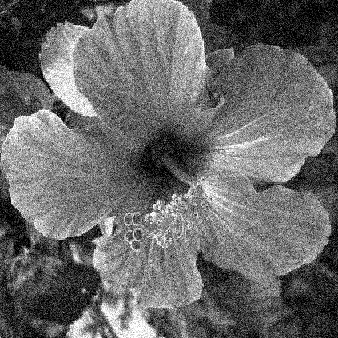}} \qquad
\subfloat[]
{\label{fig:flower_denoised_pocs}\includegraphics[width=70mm]{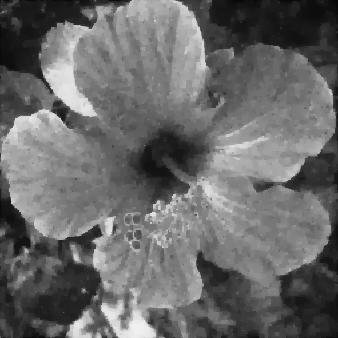}} \quad
\subfloat[]
{\label{fig:flower_denoised_chombolle}\includegraphics[width=70mm]{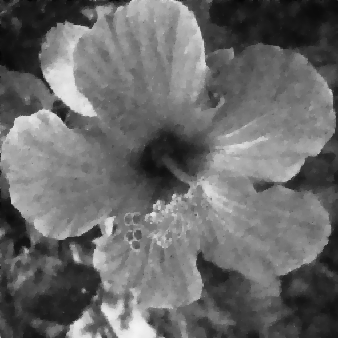}} \qquad
\caption{Sample images used in our experiments (a) Original ``Flower" image, (b) ``Flower" image corrupted with Gaussian noise with $\delta = 30$, (c) Denoised ``Flower" image, using POCS algorithm; SNR = 21.97 dB, (d) Denoised ``Flower" image, using Chambolle's algorithm; SNR = 20.89 dB.}
\label{fig:flower_example}
\end{figure}

\begin{figure}[ht]
\centering
\subfloat[]
{\label{fig:camera_org}\includegraphics[width=70mm]{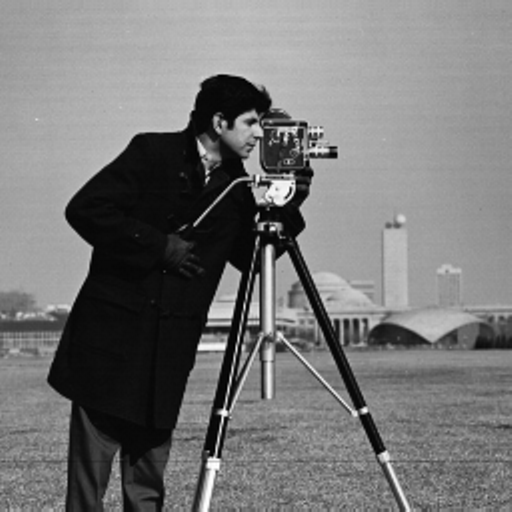}} \quad
\subfloat[]
{\label{fig:camera_noise}\includegraphics[width=70mm]{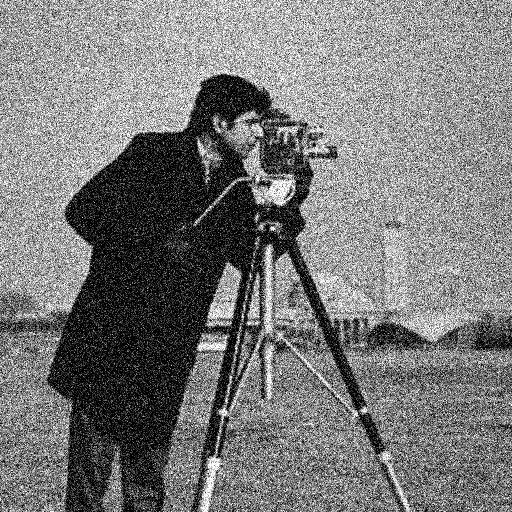}} \qquad
\subfloat[]
{\label{fig:camera_denoised_pocs}\includegraphics[width=70mm]{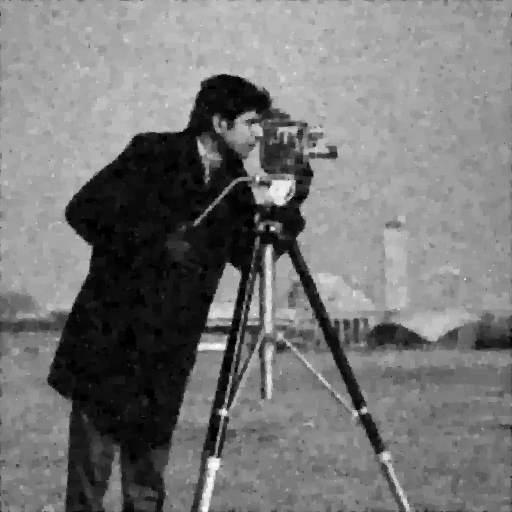}} \quad
\subfloat[]
{\label{fig:camera_denoised_chambolle}\includegraphics[width=70mm]{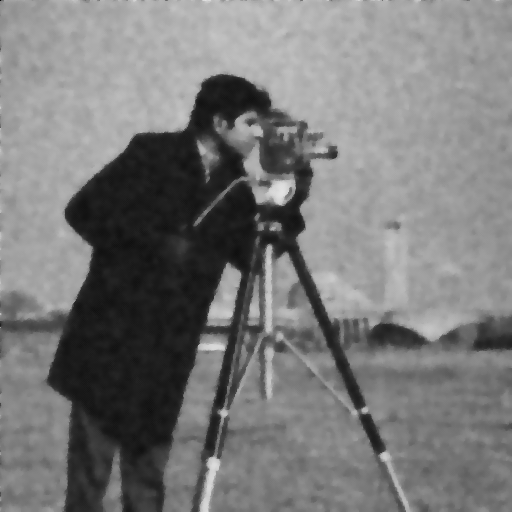}} \qquad
\caption{Sample images used in our experiments (a) Original ``Cameraman" image, (b) ``Cameraman" image corrupted with Gaussian noise with $\delta = 50$, (c) Denoised ``Cameraman" image, using POCS algorithm; SNR = 21.55 dB, (d) Denoised ``Cameraman" image, using Chambolle's algorithm; SNR = 21.22 dB.}
\label{fig:cameraman_example}
\end{figure}

Consider the ``Note" image shown in Fig. \ref{fig:note_org}. This is corrupted by a zero mean Gaussian noise with $\delta = 45$ in Fig.~\ref{fig:note_noise}. The image is restored using our method and Chombolle's algorithm \cite{Chambolle} and the denoised images are shown in Fig.~\ref{fig:note_denoised_POCS} and \ref{fig:note_denoised_chombolle}, respectively. The $\lambda$ parameter in (\ref{app:eq:cost}) is manually adjusted to get the best possible results. Our algorithm not only produce a higher SNR, but also a visually better looking image, and this is observable in two other example images in Fig. \ref{fig:cameraman_example} and \ref{fig:flower_example}, both visually and in the sense of SNR value. Solution results for other SNR levels are presented in Table \ref{tab:3}. We also corrupted ``Note" image with $\epsilon$-contaminated Gaussian noise (``salt-and-pepper noise"). De-noising results are summarized in Table \ref{tab:4}.

In Table \ref{tab:other}, de-noising results for 10 other images with two different noise levels are presented. In almost all cases our method produces higher SNR results than the de-noising results obtained using \cite{Chambolle}. The performance of the reconstruction is measured using the SNR criterion, which is defined as follows
\begin{equation}
\label{app:eq:eq6}
SNR = 20\times log_{10}(\frac{\|x\|_{2}}{\|x-x_{rec}\|_{2}}),
\end{equation}
where x is the original signal and $x_{rec}$ is the reconstructed signal.

\begin{table}[ht!]
\begin{center}
\caption{Comparison of The Results For De-noising Algorithms With Gaussian Noise For Note Image (SNRs are in dB)}
\label{tab:3}
\begin{tabular}{|c|c|c|c|}
\hline
\parbox[t]{0.8cm}{Noise \\ std }&\parbox[t]{1cm}{Input \\ SNR}&\parbox[t]{1.3cm}{\textbf{POCS}}&\parbox[t]{2.2cm}{\textbf{Chambolle}}\\\hline\hline
5&21.12&30.63&29.48\\\hline
10&15.12&25.93&24.20\\\hline
15&11.56&22.91&21.05\\\hline
20&9.06&20.93&18.90\\\hline
25&7.14&19.27&17.17\\\hline
30&5.59&17.89&15.78\\\hline
35&4.21&16.68&14.69\\\hline
40&3.07&15.90&13.70\\\hline
45&2.05&15.08&12.78\\\hline
50&1.12&14.25&12.25\\\hline
\end{tabular}
\end{center}
\end{table}

\begin{table}[ht!]
\begin{center}
\caption{Comparison of The Results For De-noising Algorithms for $\epsilon$-Contamination Noise For Note Image (SNRs are in dB)}
\label{tab:4}
\begin{tabular}{|c|c|c|c|c|c|}
\hline
$\epsilon$&$\sigma_{1}$&$\sigma_{2}$&\parbox[t]{1cm}{Input \\ SNR}&\parbox[t]{1.3cm}{\textbf{POCS}}&\parbox[t]{2.2cm}{\textbf{Chambolle}}\\\hline\hline
0.9&5&30&14.64&23.44&20.56\\\hline
0.9&5&40&12.55&21.39&17.60\\\hline
0.9&5&50&10.75&19.49&15.54\\\hline
0.9&5&60&9.29&17.61&13.82\\\hline
0.9&5&70&7.98&16.01&12.57\\\hline
0.9&5&80&6.89&14.54&11.37\\\hline\hline
0.9&10&30&12.56&22.88&19.74\\\hline
0.9&10&40&11.13&21.00&15.30\\\hline
0.9&10&50&9.85&19.35&12.47\\\hline
0.9&10&60&8.58&17.87&10.42\\\hline
0.9&10&70&7.52&16.38&8.76\\\hline
0.9&10&80&6.46&15.05&7.45\\\hline\hline
0.95&5&30&16.75&24.52&23.18\\\hline
0.95&5&40&14.98&22.59&20.44\\\hline
0.95&5&50&13.41&20.54&18.45\\\hline
0.95&5&60&12.10&18.72&16.80\\\hline
0.95&5&70&10.80&17.13&15.34\\\hline
0.95&5&80&9.76&15.63&14.11 \\\hline\hline
0.95&10&30&13.68&23.79&20.43\\\hline
0.95&10&40&12.66&22.09&15.35\\\hline
0.95&10&50&11.71&20.65&12.28\\\hline
0.95&10&60&10.72&19.10&10.22\\\hline
0.95&10&70&9.82&17.59&8.66\\\hline
0.95&10&80&8.92&16.12&7.34\\\hline
\end{tabular}
\end{center}
\end{table}

\begin{table}[ht!]
\begin{center}
\caption{Comparison of The Results For De-noising Algorithms With Gaussian Noise For Different Images With std = 30, 50 (SNRs are in dB)}
\label{tab:other}
\begin{tabular}{|c|c|c|c|c|}
\hline
\parbox[t]{0.8cm}{Images }&\parbox[t]{0.7cm}{Noise \\ std}&\parbox[t]{0.8cm}{Input \\ SNR}&\parbox[t]{1.3cm}{\textbf{POCS}}&\parbox[t]{2.2cm}{\textbf{Chambolle}}\\\hline\hline
House&30&13.85&27.43&27.13\\\hline
House&50&9.45&24.20&24.36\\\hline
Lena&30&12.95&23.63&23.54\\\hline
Lena&50&8.50&21.46&21.37\\\hline
Mandrill&30&13.04&19.98&19.64\\\hline
Mandrill&50&8.61&17.94&17.92\\\hline
Living room&30&12.65&21.21&20.88\\\hline
Living room&50&8.20&19.25&19.05\\\hline
Lake&30&13.44&22.19&21.86\\\hline
Lake&50&8.97&20.03&19.90\\\hline
Jet plane&30&15.57&26.28&25.91\\\hline
Jet plane&50&11.33&23.91&23.54\\\hline
Peppers&30&12.65&23.57&23.59\\\hline
Peppers&50&8.20&21.48&21.36\\\hline
Pirate &30&12.13&21.39&21.30\\\hline
Pirate &50&7.71&19.37&19.43\\\hline
Cameraman&30&12.97&24.13&23.67\\\hline
Cameraman&50&8.55&21.55&21.22\\\hline
Flower&30&11.84&21.97&20.89\\\hline
Flower&50&7.42&19.00&18.88\\\hline\hline
Average&30&13.11&23.18&22.84\\\hline
Average&50&8.69&20.82&20.70\\\hline
\end{tabular}
\end{center}
\end{table}

\begin{figure}[ht]
\centering
\subfloat[House]
{\label{fig:house}\includegraphics[width=40mm]{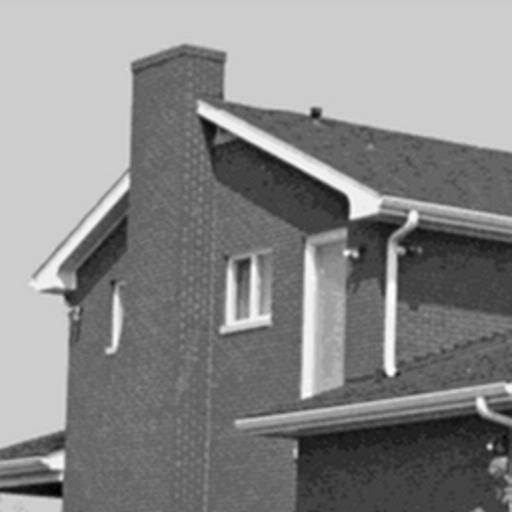}} \quad
\subfloat[Jet plane]
{\label{fig:jetplane}\includegraphics[width=40mm]{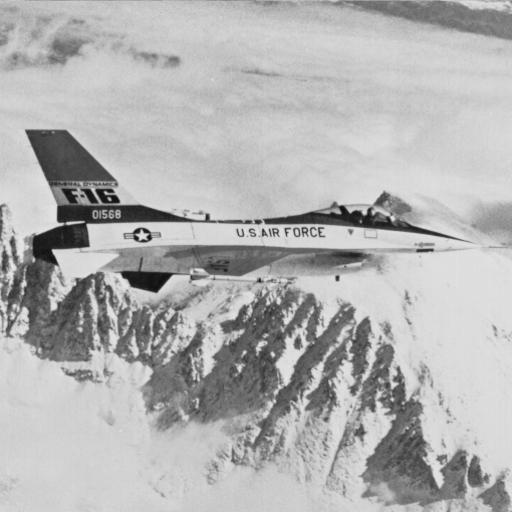}} \quad
\subfloat[Lake]
{\label{fig:lake}\includegraphics[width=40mm]{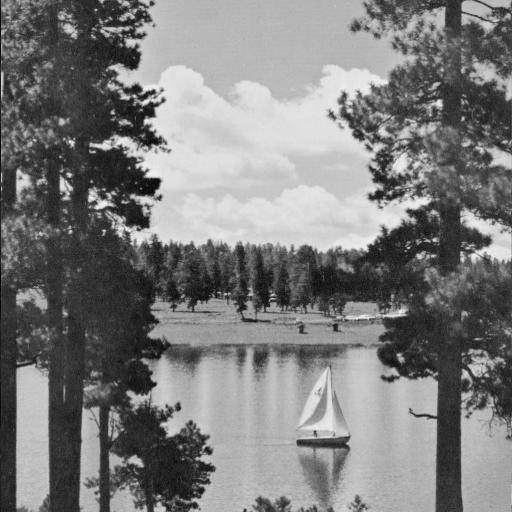}} \quad
\subfloat[Lena]
{\label{fig:lena}\includegraphics[width=40mm]{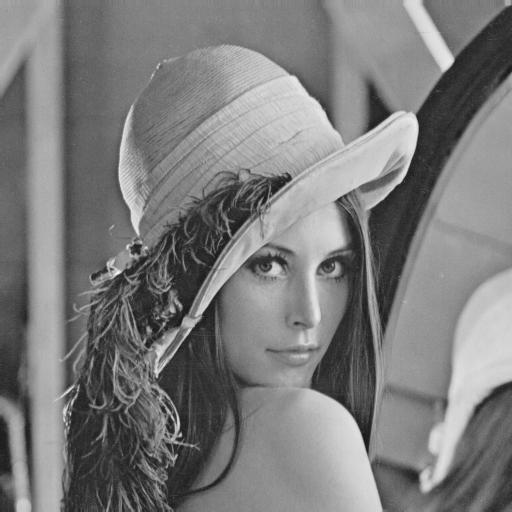}} \quad
\subfloat[Living room]
{\label{fig:livingroom}\includegraphics[width=40mm]{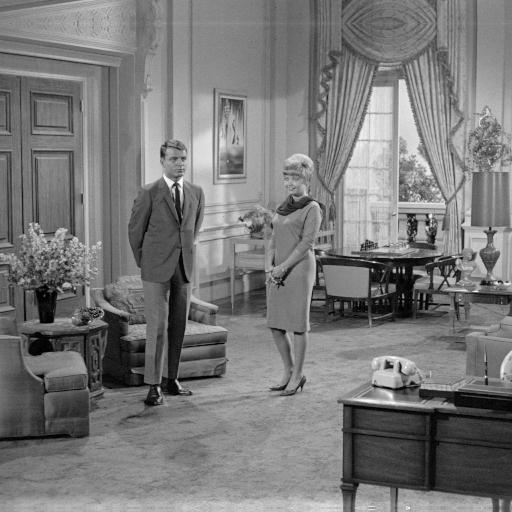}} \quad
\subfloat[Mandrill]
{\label{fig:mandril}\includegraphics[width=40mm]{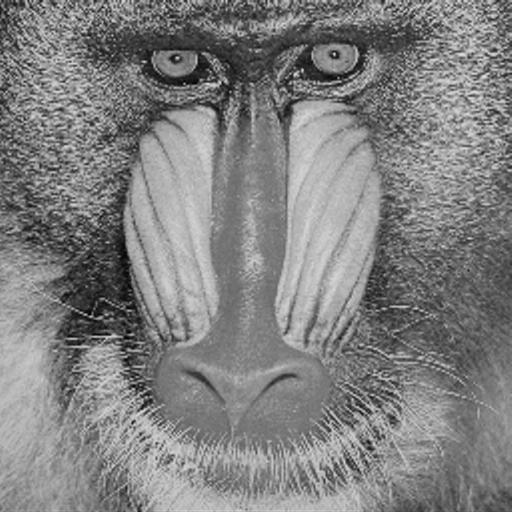}} \quad
\subfloat[Peppers]
{\label{fig:peppers}\includegraphics[width=40mm]{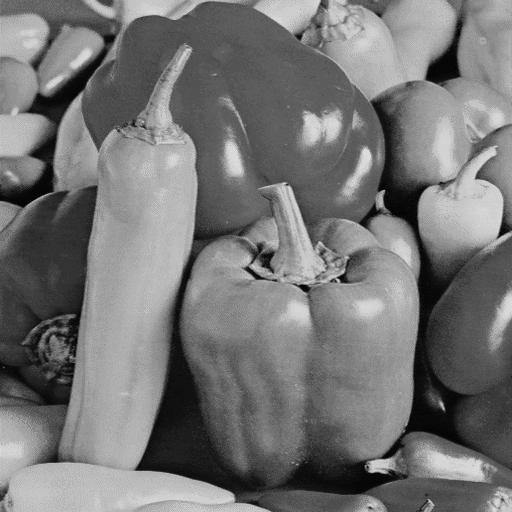}} \quad
\subfloat[Pirate]
{\label{fig:pirate}\includegraphics[width=40mm]{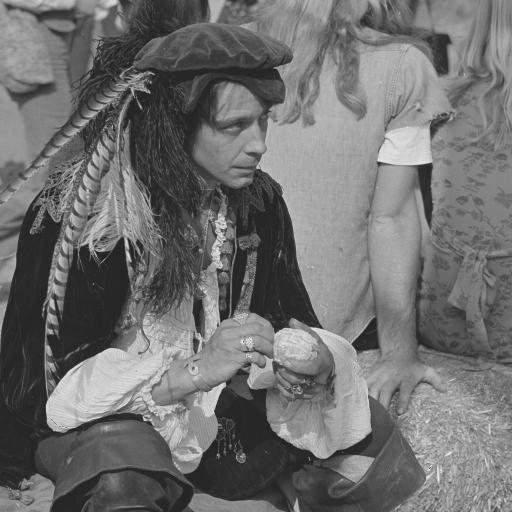}}
\caption{Sample images used in our experiments (a) House, (b) Jet plane, (c) Lake, (d) Lena, (e) Living room, (f) Mandrill, (g) Peppers, (h) Pirate.}
\label{fig:otherOrg}
\end{figure}

\begin{figure}[ht!]
\begin{center}
\noindent
\includegraphics[width=150mm]{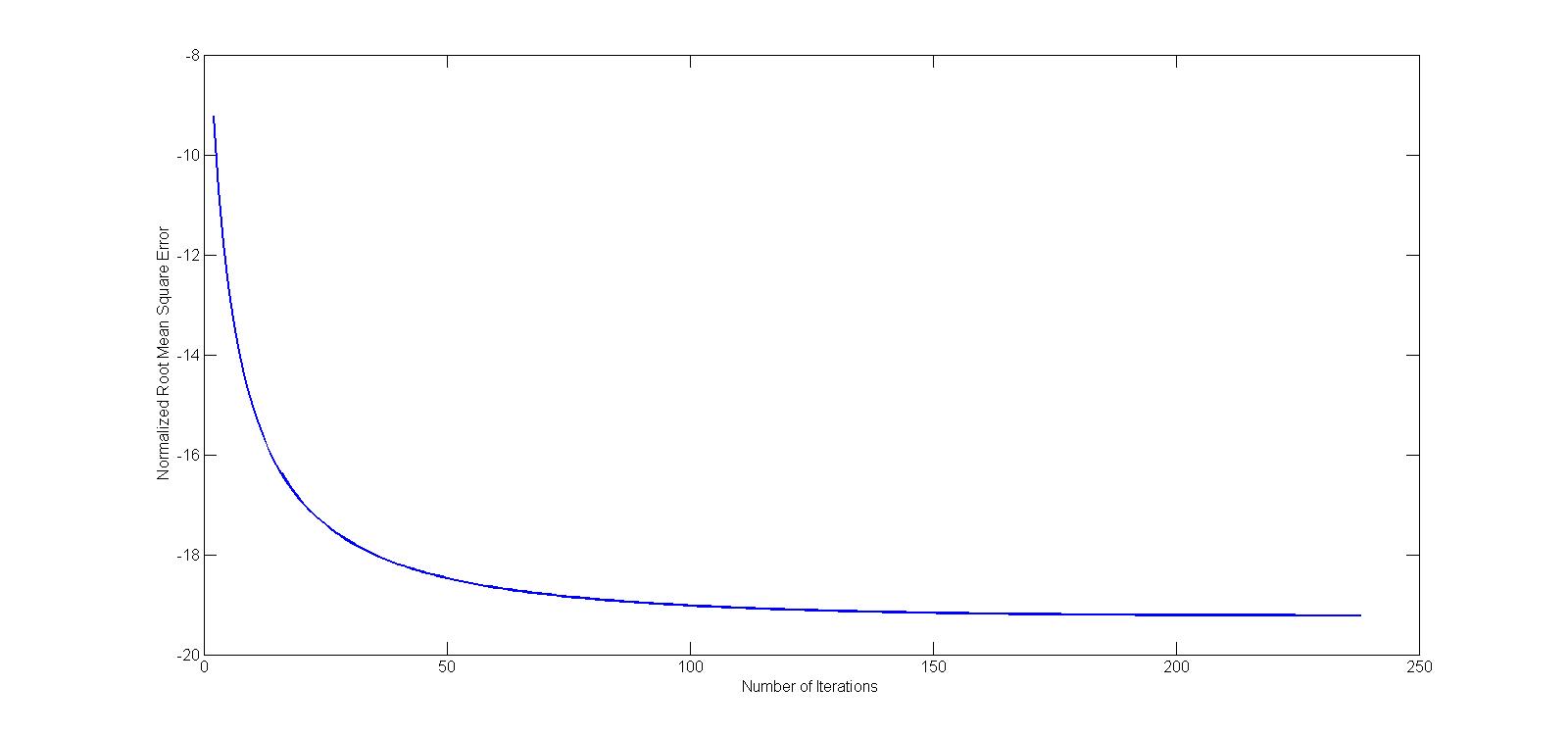}
\caption{Normalized Root Mean Square Error in each iteration.}
\label{app:NRMSE}
\end{center}
\end{figure}

Moreover, to illustrate the convergence process of the proposed algorithm with ``Note" image corrupted with Gaussian noise with standard deviation equal to 25, we calculated the Normalized Root Mean Square Error ($NRMSE$) as:
\begin{equation}
\label{app:eq:eq8}
NRMSE(i) = \frac{\|x_{i} - x\|}{\|x\|},
\end{equation}
where $x_{i}$ is the de-noised image in $i^{th}$ step, and $x$ is the original image.As in Fig. $\ref{app:NRMSE}$, the NRMSE is decreasing from -9 dB to -19 dB.

\begin{figure}[ht!]
\begin{center}
\noindent
\includegraphics[width=150mm]{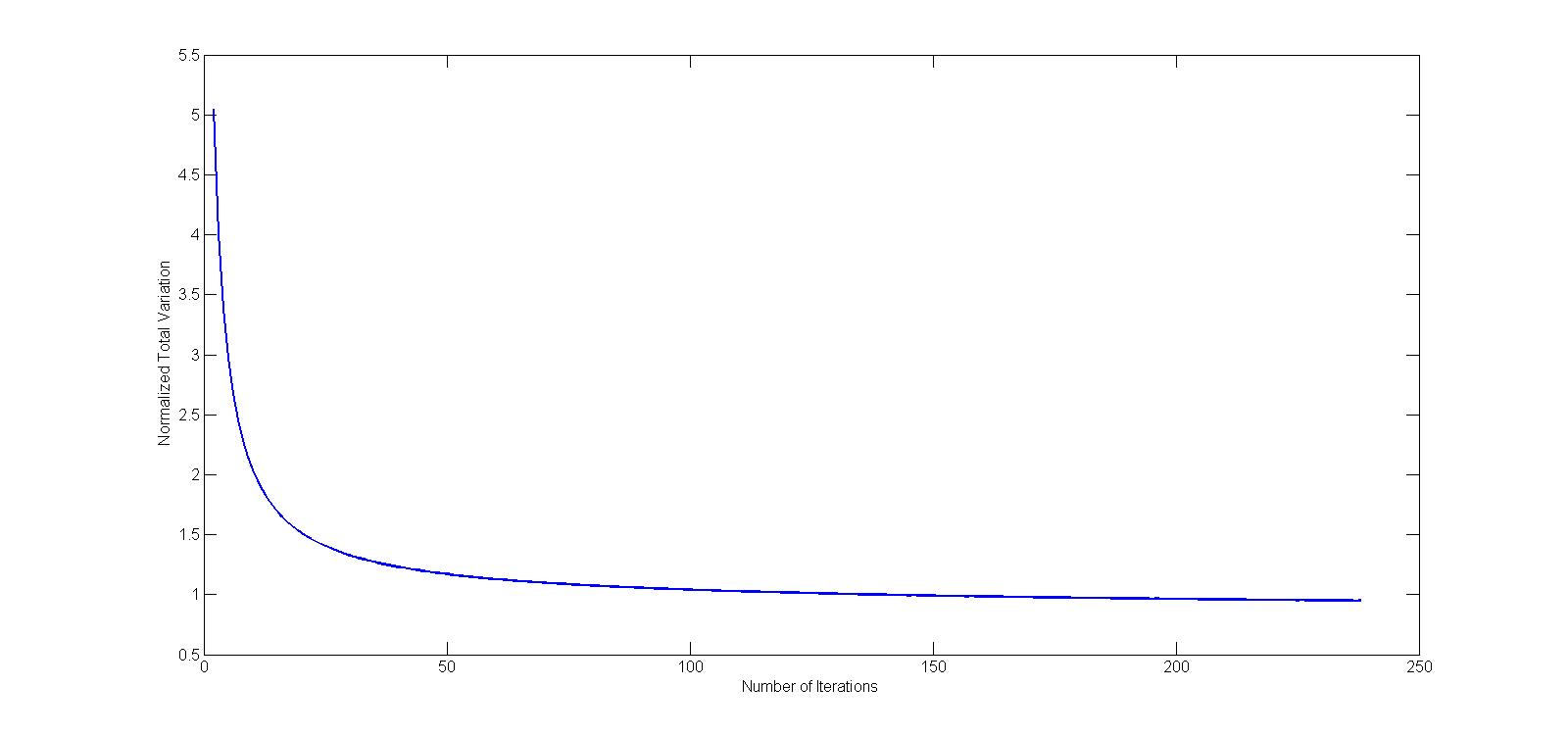}
\caption{Normalized Total Variation in each iteration.}
\label{app:NTV}
\end{center}
\end{figure}

\begin{figure}[ht!]
\begin{center}
\noindent
\includegraphics[width=150mm]{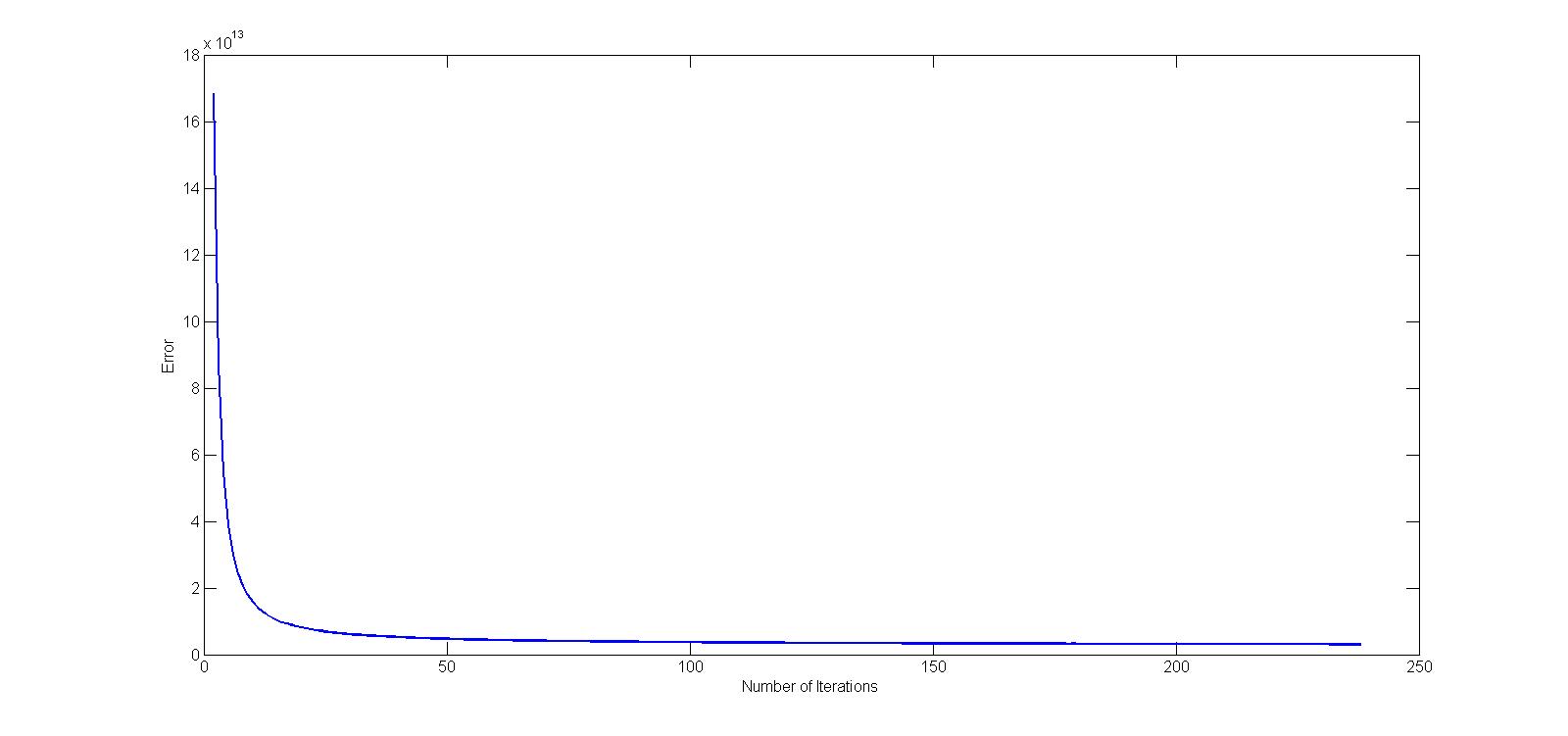}
\caption{Convergence error curve in each iteration.}
\label{app:ErrorCurve}
\end{center}
\end{figure}

As another experiment to show the convergence behavior, we display the Normalized Total Variation (NTV) as:
\begin{equation}
\label{app:eq:eq9}
NTV(i) = \frac{TV(x_{i})}{TV(x)},
\end{equation}
where $x_i$ and $x$ are the same as the $NRMSE$. As is obvious in Fig. $\ref{app:NTV}$, the Normalized TV curve has converged to almost 1, which is demonstrating the successful convergence. As the last measurement of convergence, consider the error value in each step of the iteration; Fig. $\ref{app:ErrorCurve}$, shows the error value in each step. These three curves shows that the iterations converge to the desired solution roughly around the $100_{th}$ iteration. The convergence iteration number depends on the noise level, as here for $\delta = 25$ its almost 100.

\subsection{Compressive sensing}
\label{subsec:Simulation Compressive sensing}
\begin{figure}[ht!]
\begin{center}
\noindent
\includegraphics[width=150mm]{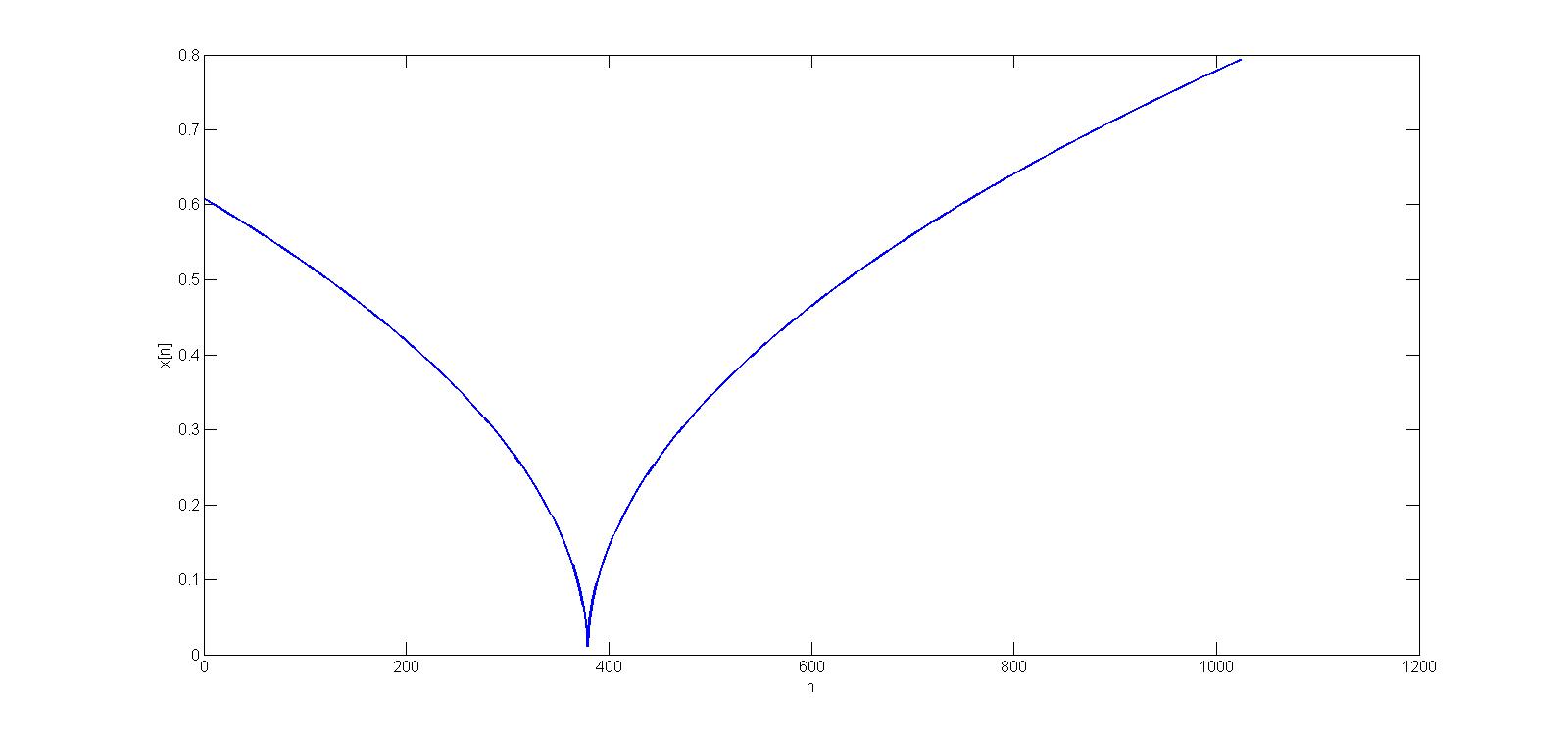}
\caption{The cusp signal with N = 1024 samples.}
\label{app:cusp}
\end{center}
\end{figure}

For the validation and testing of the proposed algorithm, experiments are carried out with one-dimensional (1D) signals, and two-dimensional (2-D) signals, including 30 different images. For 1-D signal, we've done experiments with the cusp signal, which consists of 1024 samples, and is shown in Figure \ref{app:cusp}. In the DCT domain, the cusp signal can be approximated sparsely . The sparse random signals consisting of 4 and 25 randomly located non-zero samples with random values, are composed of 128 and 256 samples, respectively. In all the experiments, the measurement matrices $\phi$ are chosen as Gaussian random matrices.

\begin{figure}[ht]
\centering
\subfloat[]
{\label{fig:cusp204}\includegraphics[width=120mm]{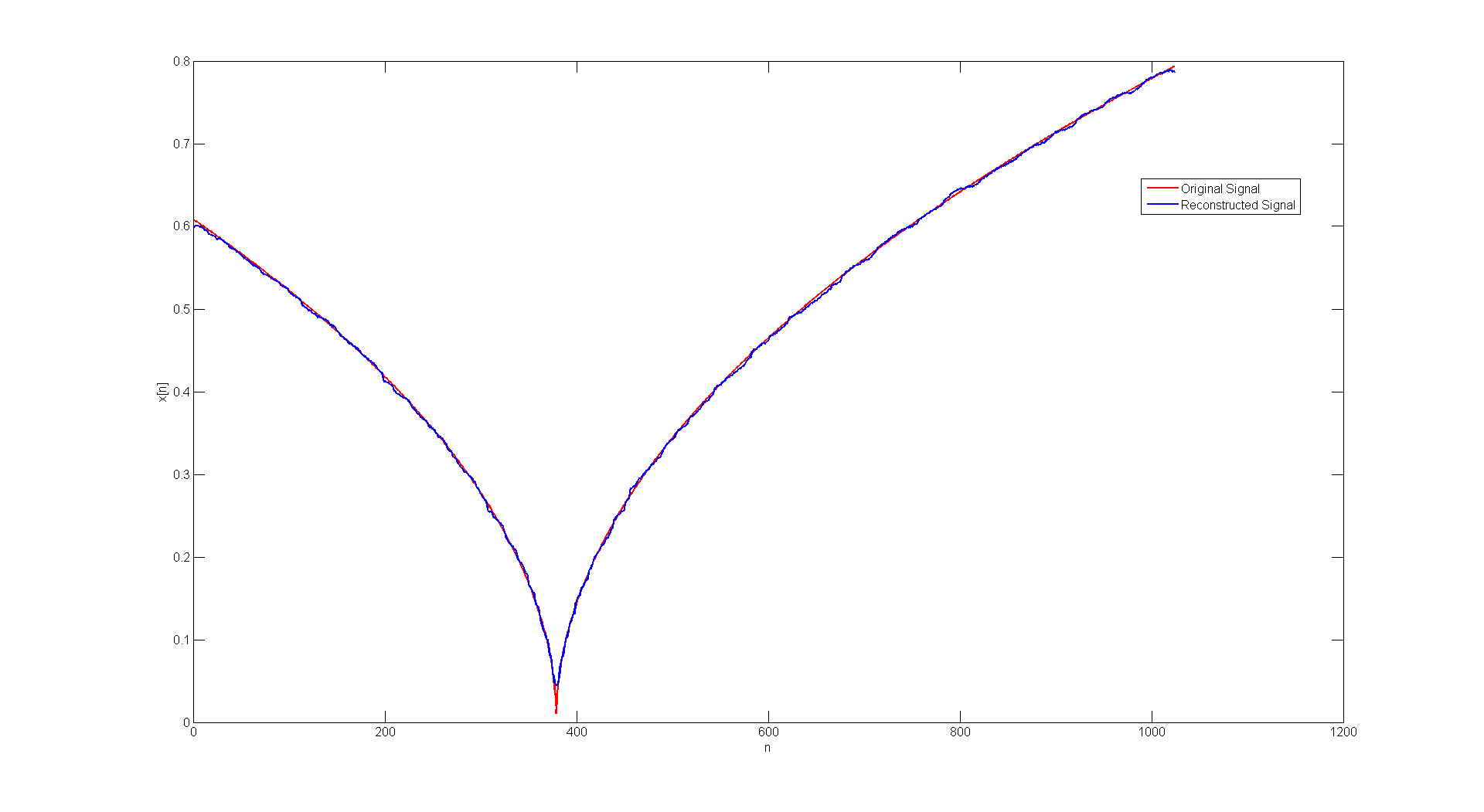}} \quad
\subfloat[]
{\label{fig:cusp717}\includegraphics[width=120mm]{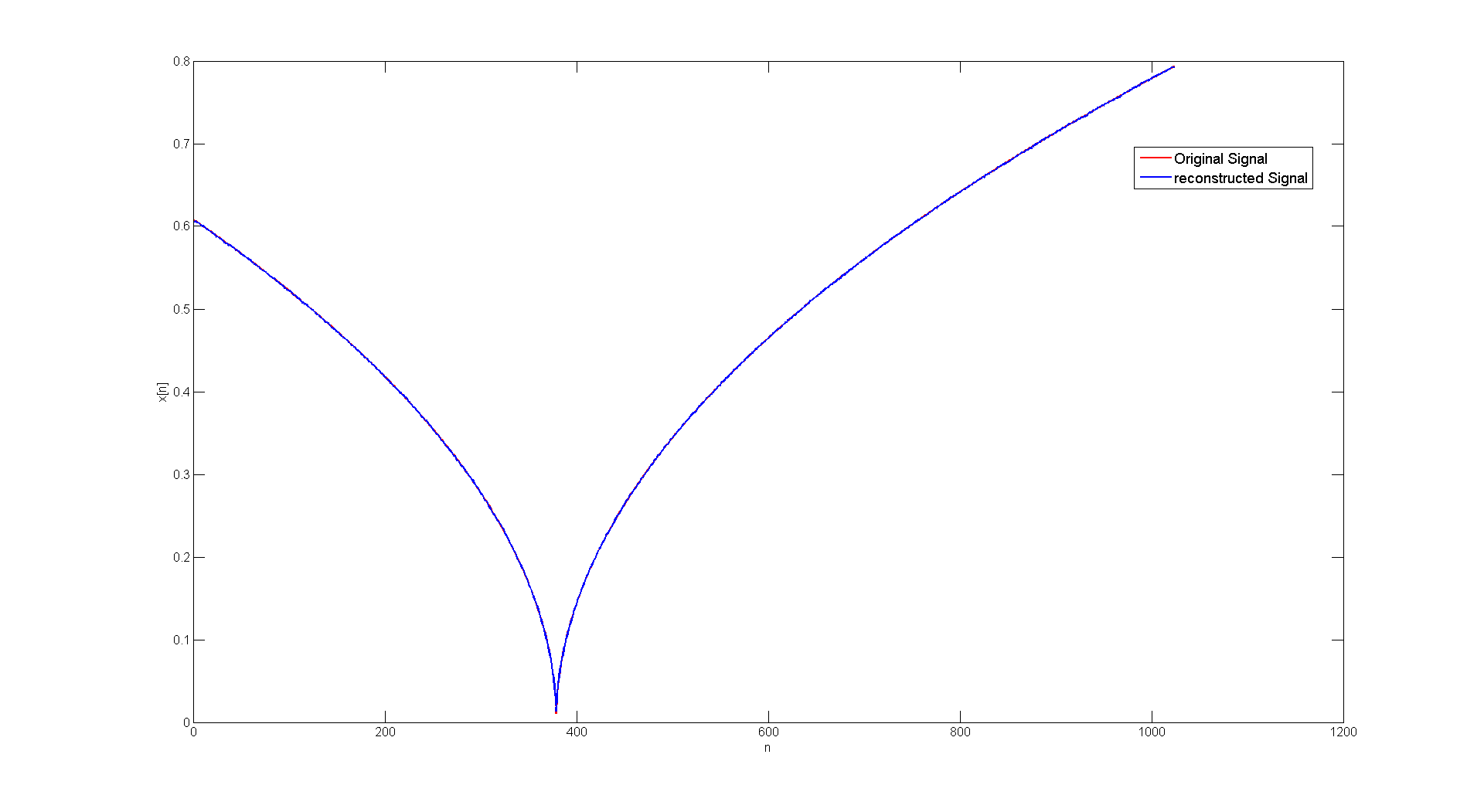}}
\caption{The reconstructed cusp signal with N = 1024 samples, for (a) 204 measurements (SNR = 45), and (b) 717 measurements (SNR = 58).}
\label{fig:otherOrg}
\end{figure}

As the first set of experiments, the original cusp signal is reconstructed, with M = 204, 717 measurements and, when M = 24, 40 measurements are taken from the S = 5 random signal with 128 samples. The reconstructed signals using the TV cost functional based algorithm are shown in Figures $\ref{fig:cusp204}$, and $\ref{fig:cusp717}$, 6(a), and 6(b). In DCT domain, the cusp signal has 76 coefficients with large magnitude. Therefore, it can be approximated as S = 76 sparse signal in the DCT domain.

After reconstructing the original cusp signal using the proposed algorithm with M = 204, and 717 measurements, we've obtained results of 45 and 58 dB, SNR values. In the case of the experiment with random signals, the proposed method reconstructed the original signal with a slight error when using 30 measurement; however, the samples are detected appropriately. Moreover, it is perfectly reconstructed the original signal when using 50 measurements.

\begin{figure}[ht!]
\begin{center}
\noindent
\includegraphics[width=150mm]{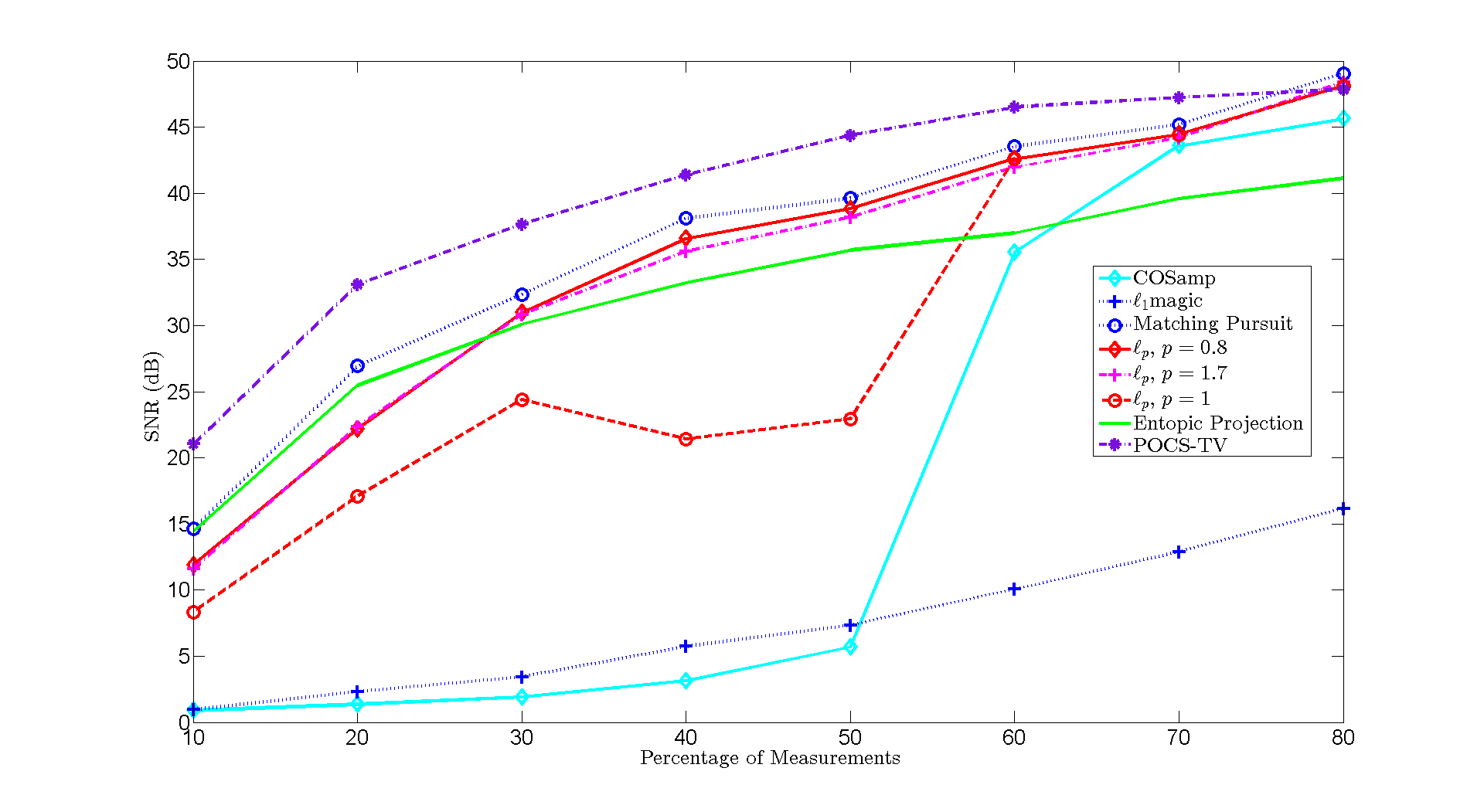}
\caption{The reconstructed cusp signal with N = 256 samples.}
\label{app:cusp_recon}
\end{center}
\end{figure}

In the next set of experiments, we implemented the proposed algorithm in 2-dimension (2D) and applied it to six well known images from the image processing literature and 24 images from the ``Kodak True Color images" database \cite{toolbox}. Table \ref{tab:comp} represents the SNR values for compressive sensing with two algorithms for these images with two block-size of $32\times32$ and $64\times64$ for both algorithms. The images in Kodak dataset are 24 bit per pixel color images. We first transformed all the color images into YUV color space and used the 8 bit per pixel luminance component (Y channel) of the images in our tests. Figure \ref{app:cusp_recon} representing the SNR values for some compressive sensing algorithms with eight different percentages of measurements, ranging from 10\% up to 80\%.

We compared our results with the block based compressed sensing algorithm given in \cite{fowler}. Therefore, we divided the image into blocks and reconstructed those blocks individually. Random measurements, which are 30\% of the total number of points in images, are used in tests on both the proposed algorithm and Fowler $et. al.$’s method. On average, for $64\times64$, and $32\times32$ blocks, we achieved approximately 1.24 dB, and 0.42 dB higher SNR respectively, compared to the algorithm given in \cite{fowler}, as shown in Table \ref{tab:comp}.

\begin{table}[ht!]
\begin{center}
\caption{Comparison of The Results For Compressive Sensing With Fowlers algorithm and Our Algorithm}
\label{tab:comp}
\begin{tabular}{|c|c|c|c|c|}
\hline
\parbox[t]{0.8cm}{Image}&\parbox[t]{1cm}{Fowler \\ B=32}&\parbox[t]{1cm}{Fowler \\ B=64}&\parbox[t]{1.5cm}{POCS \\ B=32}&\parbox[t]{1.5cm}{POCS \\ B=64}\\\hline\hline
Kodak(ave.)&21.40&21.36&21.96&22.74\\\hline
Mandrill&16.47&16.77&16.65&16.96\\\hline
Lena&26.82&26.71&26.02&26.86\\\hline
Barbara&20.05&20.40&18.21&18.62\\\hline
Peppers&24.66&24.46&27.06&27.93\\\hline
Goldhill&22.78&23.44&23.64&24.24\\\hline
Average&21.60&21.53&22.02&22.77\\\hline
\end{tabular}
\end{center}
\end{table}

In the last set of experiments, we compared our reconstruction results with 4 well known CS reconstruction algorithms from the literature; CoSamp \cite{Needell2009301}, $\ell_1$magic \cite{Candes1}, Matching Pursuit (MP) \cite{Mallat}, and $\ell_p$ optimization based CS reconstruction \cite{Chartrand2} algorithms. In comparison to the $\ell_p$ optimization based CS reconstruction algorithm, we used three different values for p; p = [0.8, 1, 1.7]. With $p = 1$, the algorithm solves the problem given in \ref{app:eq:eq5}, which is the $\ell_1$ norm optimization problem.

\begin{figure}[ht!]
\begin{center}
\noindent
\includegraphics[width=150mm]{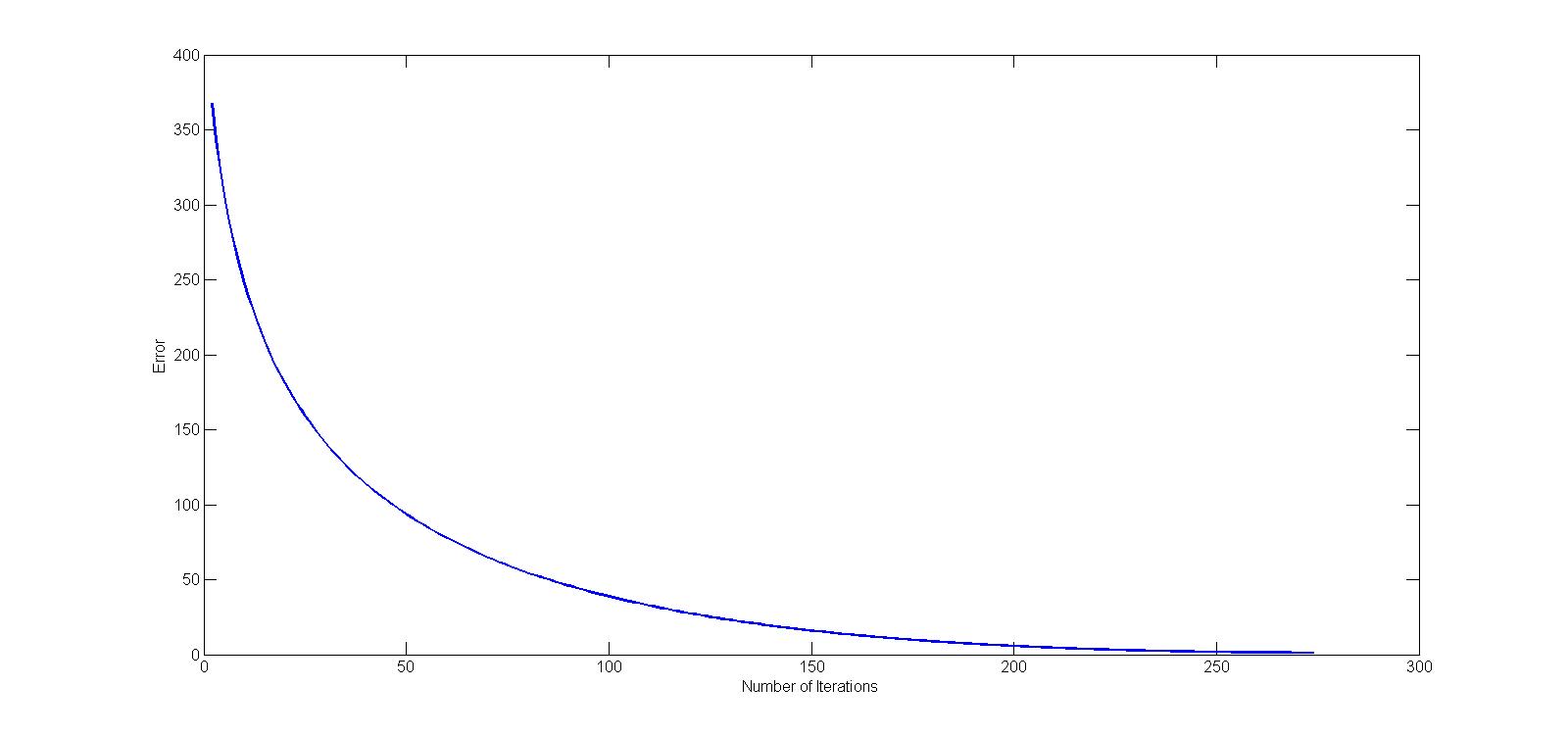}
\caption{The convergence curve for cusp signal with N = 1024 samples, and M = 307 measurements.}
\label{app:cusp_recon_307}
\end{center}
\end{figure}

In this set of experiments, to reconstruct the original signal, we've also implemented the both algorithms with different number of measurements as in previous set of experiments ranging from 10\% to 80\% of the total number of the samples of the 1D signal. Then, the SNR values are measured between the original signals and the reconstructed ones. The main region of interest in these experiments is 20\% − 60\% range. The results of the tests with cusp signal are presented in Fig. \ref{app:cusp_recon}. The proposed algorithm in almost all of the cases performed better than other algorithms. As another measurement of convergence, we calculated the error in each step of the iterations, and as in Fig. $\ref{app:cusp_recon_307}$, the error cure converges almost to zero for $M = 307$ (30\% of signal length). It is important to note that, the cusp signal is not sparse; however, since the coefficients in most of the transform domains are not zero, but close to zero, then it is compressible \cite{cevher3}. Therefore, the sparsity level of the test signals are not known exactly beforehand.

\begin{figure}[ht]
\centering
\subfloat[]
{\label{fig:house}\includegraphics[width=70mm]{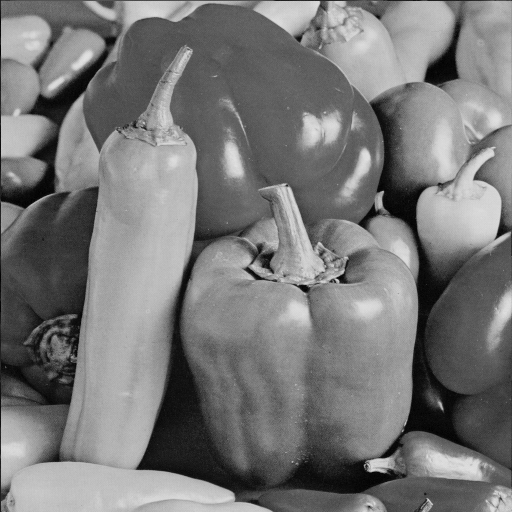}} \quad
\subfloat[]
{\label{fig:jetplane}\includegraphics[width=70mm]{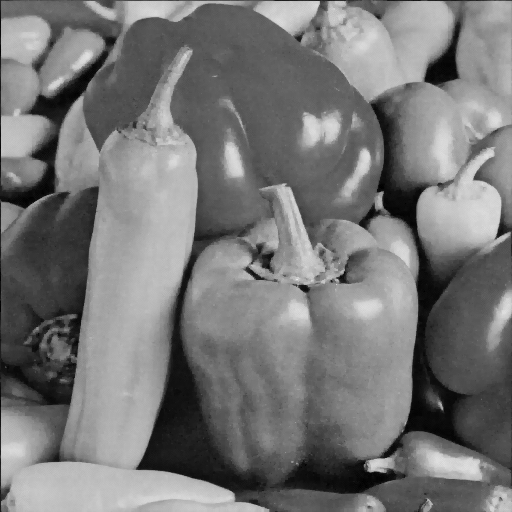}} \qquad
\subfloat[]
{\label{fig:lake}\includegraphics[width=70mm]{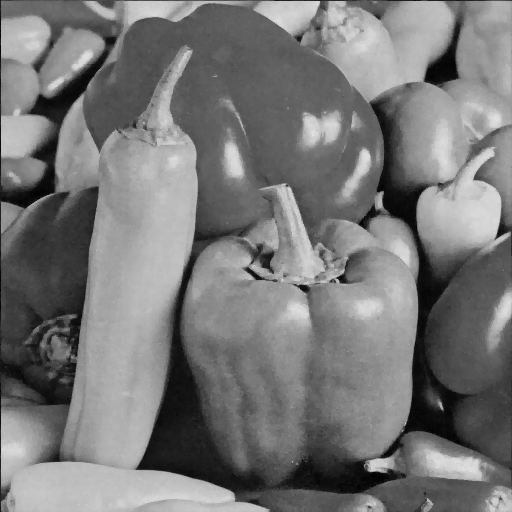}}
\caption{Sample results for peppers image with our algorithm (a) Original Image , (b) Reconstructed image with 32x32 blocks, SNR = 27.06, (c) Reconstructed image with 64x64 blocks, SNR = 27.93.}
\label{fig:otherOrg}
\end{figure}

\begin{figure}[ht]
\centering
\subfloat[]
{\label{fig:house}\includegraphics[width=70mm]{peppers.jpg}} \quad
\subfloat[]
{\label{fig:jetplane}\includegraphics[width=70mm]{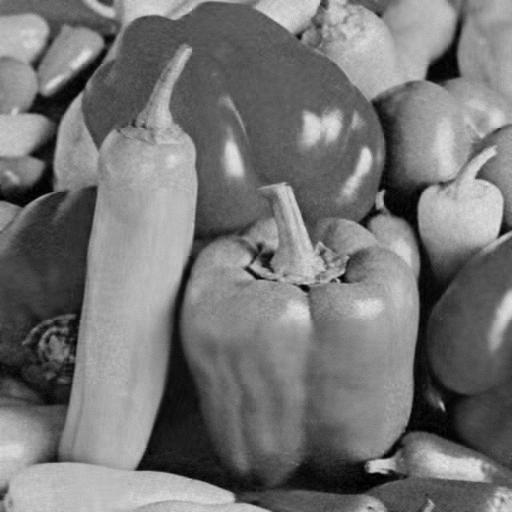}} \qquad
\subfloat[]
{\label{fig:lake}\includegraphics[width=70mm]{BCS_peppers_64.jpg}}
\caption{Sample results for peppers image with Fowler's algorithm (a) Original Image , (b) Reconstructed image with 32x32 blocks, SNR = 24.66, (c) Reconstructed image with 64x64 blocks, SNR = 24.46.}
\label{fig:otherOrg}
\end{figure}

\section{Conclusion}
\label{sec:Conclusion}
A new de-noising method based on the epigraph of the TV function is developed. The solution is obtained using POCS. The new algorithm does not need the optimization of the regularization parameter.

\clearpage
\bibliographystyle{IEEEtran}
\bibliography{PhdReferences}

\end{document}